%% file: schrittesser-on-horowitz-shelah.tex
\documentclass[a4paper,10pt,reqno]{amsart}

\usepackage[applemac]{inputenc}
\usepackage[T1]{fontenc}
\usepackage[english]{babel}
\usepackage{csquotes}

\usepackage{amssymb}
\usepackage{amsmath,amsthm}

\usepackage[shortlabels]{enumitem}

\usepackage[pdfauthor={David Schrittesser},
            pdftitle={},
            pdfproducer={Latex with hyperref}]{hyperref}

\include{preamble}

\newcommand{\order}[1][f]{\mathbin{\prec^{#1}}}

\DeclareMathOperator{\code}{\#}

\DeclareMathOperator{\powerset}{\mathcal{P}}

\DeclareMathOperator{\domain}{B}

\DeclareMathOperator{\ed}{e}
\DeclareMathOperator{\edb}{\ddot e}
\DeclareMathOperator{\edh}{\dot e}
\DeclareMathOperator{\ho}{C}

\DeclareMathOperator{\predT}{\mathcal T}

\title{On Horowitz and Shelah's Borel maximal eventually different family}

\author{David Schrittesser}
\address{Department of Mathematical Sciences, University of Copenhagen, Universitetsparken 5, 2100 Copenhagen \O, Denmark}
\email{david@logic.univie.ac.at}



\begin{document}

\newcommand{\linspan}[1]{\hull( #1 )}

\begin{abstract}
 We show there is a closed (in fact effectively closed, i.e., $\Pi^0_1$) eventually different family (working in ZF or less).
\end{abstract}

\maketitle

\section{Introduction}

\paragraph{A} 
We call a set $\mathcal E$ an \emph{eventually different family} (of functions from $\nat$ to $\nat$) if and only if 
$\mathcal E \subseteq {}^\nat\nat$ and
any two distinct $f_0, f_1 \in \mathcal E$ are \emph{eventually different}, i.e., $\{ n \setdef f_0(n) = f_1(n) \}$ is finite;
such a family is called \emph{maximal} if and only if 
it is maximal with respect to inclusion among eventually different families (we abbreviate \emph{maximal eventually different family} by \emph{medf}).

 \medskip
 
In \cite{medf-borel} Horowitz and Shelah prove the following (working in ZF).
\begin{thm}[\cite{medf-borel}]\label{t.borel} 
There is a $\Delta^1_1$ (i.e., effectively Borel) maximal eventually different family.
\end{thm}
This was surprising as the analogous statement is false in many seemingly similar situations: e.g., infinite so-called \emph{mad} families cannot be analytic \cite{mathias.1967} (see also \cite{tornquist.2015}). In a more recent, related result \cite{mcg-borel} Horowitz and Shelah obtain a $\Delta^1_1$ \emph{maximal cofinitary group}.

\medskip

In this note we present a short and elementary proof of the following improvement of Theorem~\ref{t.borel}:
\begin{thm}\label{t.main}
 There is a $\Pi^0_1$ (i.e., effectively closed) maximal eventually different family.
 \end{thm}
To prove this we first define an \emph{medf} in a simpler manner than \cite{medf-borel} (its defining formula will be $\Sigma^0_3\vee \Pi^0_3$). We then show that we can produce from any arithmetic \emph{medf} a new \emph{medf} whose definition contains one less existential quantifier. The main result follows.

\medskip

{\it Note:} 
Theorem~\ref{t.main} was announced by Horowitz and Shelah without proof in \cite{medf-borel}; the proof in the present paper was found by the author while studying their construction of a $\Delta^1_1$ \emph{medf} in \cite{medf-borel}.

In a related paper \cite{medf-bounded} the present author presents a further simplification of the construction and positively answers the following question of Asger Törnquist \cite{personal}: Given $F\colon\nat \to \nat$ such that $\liminf_{n\to\infty} F(n) = \infty$, is there a Borel (or even compact) \emph{medf} in the \emph{restricted} space 
$\mathcal N_F = \{g \in {}^\nat\nat \setdef (\forall n\in \nat)\; g(n) \leq F(n) \}$?

\medskip

\paragraph{B} We fix some notation and terminology (generally, our reference for notation is \cite{kechris}). 
 `$\exists^\infty$' means `there are infinitely many\dots,'
${}^\nat \nat$ means the set of functions from $\nat$ to $\nat$ and ${}^{<\nat} \nat$ means the set of finite sequences from $\nat$; we write $\lh(s)$ for the length of $s$ when $s \in {}^{<\nat}\nat$.
For  $s,t\in {}^{n}\nat$, $s\conc t$ is the \emph{concatenation} of $s$ and $t$, i.e., the unique $u \in {}^{\lh(s)+\lh(t)}\nat$
such that $s \subseteq u$ and $(\forall k<\lh(t)) \; u(\lh(s)+k) = t(k)$.

We write $f_0 \mathbin{=^\infty} f_1$ to mean that $f_0$ and $f_1$ are \emph{not} eventually different (they are infinitely equal).
 Two sets $A, B \subseteq \nat$ are called \emph{almost disjoint} if and only if $A \cap B$ is finite, and
an \emph{almost disjoint family} is a set $\mathcal A \subseteq \powerset(\nat)$ any two elements of which are almost disjoint. 

Qualifications like ``\dots is recursive (i.e., computable) in\dots'' are applied to subsets of $\Hhier(\omega)$, the set of hereditarily finite sets.
Consult \cite{moschovakis,mansfield,kechris} for more on the (effective) Borel and projective hierarchies, i.e., on $\Pi^0_1$, $\Pi^0_1(F)$,  $\Delta^1_1$, \dots{} sets.

All results in this paper can be derived in ZF (or in fact, in a not so strong subsystem of second order arithmetic).

\medskip

\paragraph{C}
This note is organized as follows.
In Section~\ref{s.recipe} we make some motivating observations, leading to Lemma~\ref{l.abstract}  which gives an abstract recipe for creating maximal eventually different families.
We take the opportunity to give a rough sketch of the proof of Theorem~\ref{t.borel} as given by Horowitz and Shelah in \cite{medf-borel}.

We then give a simpler construction instantiating the recipe from Lemma~\ref{l.abstract} and 
yielding a \emph{medf} which is $\Sigma^0_3\vee\Pi^0_3$ in Section~\ref{s.construction}.
Lastly, we show how to get rid of all existential quantifiers in Section~\ref{s.mangling}. This requires mangling the family, but the new family is still  maximal eventually different.

\medskip

{\it Acknowledgements: The author gratefully acknowledges
 the generous support from the DNRF Niels Bohr Professorship of Lars Hesselholt. 
}

\section{The recipe}\label{s.recipe}

\begin{dfn}\label{d.ed}
Fix a computable (i.e., $\Delta^0_1$) bijection $ n \mapsto s_n$ of $\nat$ with ${}^{<\nat}\nat$ and write $s \mapsto \code s$ for its inverse.
Given $f \colon \nat \to \nat$, let $e(f) \colon \nat \to \nat$ be the function defined by
\[
\ed(f) (n) = \code f\res n.
\]
\end{dfn}
Clearly $\{ \ed(f) \setdef f\in {}^{\nat}\nat \}$ is an eventually different family. 
At first sight, it may seem a naive strategy to make it also maximal by varying the definition of $\ed(f)$ so that it leaves $f$ intact on some infinite set. But this is just how \cite{medf-borel} succeeds.
\begin{dfn}\label{d.edb}
Let $f\colon \nat \to \nat$.
\begin{enumerate}[ref=\Alph*,label=\Alph*.]
\item\label{d.edb.domain} Let $\domain(f) = \{ 2n +1 \setdef s_n\subseteq f \}$.
\item For a set $B\subseteq \nat$, let $\edb(f,B)\colon \nat \to \nat$ be the function defined by
\[
\edb (f,B) (n) = \begin{cases} 
 f(n) & \text{if $n \in B$,}\\
\code f\res n & \text{if $n \notin B$.}
\end{cases}
\]
\end{enumerate}
\end{dfn}
\begin{rem}\label{r.f.recursive}
Note for later that $f$ is recursive in $\edb(f,\domain(f))$ as  $\edb(f,\domain(f))\res 2\nat = \ed(f)\res 2\nat$.
\end{rem}

The family $\mathcal E_0= \{ \edb(f,\domain(f)) \setdef f\in {}^{\nat}\nat \}$ is \emph{spanning}, i.e., $(\forall h \in {}^\nat\nat)(\exists  g \in \mathcal F) \; h \mathbin{=^\infty} g$.
Interestingly, $\mathcal E_0$ is also in some sense close to being eventually different:
For if 
$\edb(f,\domain(f))(n) = \edb(f',\domain(f'))(n)$ for infinitely many $n$, almost all of these $n$ must lie in $\domain (f) \cup \domain(f')$
and hence as $\{ \domain(f) \setdef f\in {}^{\nat}\nat \}$ is an almost disjoint family,
\[
(\exists^\infty n \in \domain(f))\; f(n) = \ed(f') (n)
\]
or the same holds with $f$ and $f'$ switched.

 \medskip
 
The brilliant idea of Horowitz and Shelah is the following:
Ensure maximality with respect to $f$ which look like $\ed(f')$ on an infinite set using $\ed(f')$;
 restrict the use of $\edb$ to $f$ which \emph{don't} look like they arise from $\ed$ on some infinite subset of $\domain(f)$ to avoid the situation described above. 
We make these ideas precise in the following definition and in Lemma~\ref{l.abstract} below.

\begin{dfn}
Let a function $f \colon \nat \to \nat$ and $X \subseteq \nat$ be given.
We say \emph{$f$ is $\infty$-coherent on $X$} if and only if there is $f' \in {}^\nat\nat$ and infinite $X' \subseteq X$ such that $f\res X' = e(f')\res X'$.
\end{dfn}
We can now give a general recipe for constructing a \emph{medf}.
\begin{lem}\label{l.abstract}
Suppose that $\predT \subseteq {}^\nat \nat$ and  $\ho\colon {}^\nat\nat \to \powerset(\nat)$ is a function such that
\begin{enumerate}[label=(\Alph*),ref=\Alph*]
\item\label{T.2}  If $f \notin \predT$, there is an infinite set $X' \subseteq \ho(f)$ and $f'\in {}^\nat \nat$ such that $f \res X' = \ed(f') \res X'$; i.e.,  $f$ is $\infty$-coherent on $\ho(f)$.
\item\label{T.1}  If $f \in \predT$,   for no $f' \in {}^\nat\nat$ does $f$ agree with $e(f')$ on infinitely many points in $\ho(f)$; i.e., $f$ is not $\infty$-coherent on $\ho(f)$.
\item\label{adf} $\{ \ho(f) \setdef f  \in \predT \}$ is an almost disjoint family.
\end{enumerate}
Then
\[
\mathcal E = \{ \edb(f,\ho(f)) \setdef f  \in \predT \} \cup \{ \ed(f) \setdef f  \notin \predT \}
\]
is a maximal eventually different family.
\end{lem}
Of course the challenge here is to define $\ho$ and $\predT$ so that $\mathcal E$ is $\Delta^1_1$; before we discuss this aspect, we prove the lemma.

For the sake of this proof it will be convenient to define the map $\edh\colon {}^\nat\nat\to{}^\nat\nat$ as follows:
For $f \in {}^\nat\nat$ let $\edh(f)$ be the function defined by
\begin{equation}\label{e.e.ddot}
\edh(f) = \begin{cases}
\edb(f,\ho(f)) & \text{ if $f \in \predT$,}\\
\ed(f)             & \text{ otherwise.}
\end{cases}
\end{equation}
Clearly $\mathcal E =\{ \edh(f) \setdef f \in {}^\nat\nat\}$.
\begin{proof}[Proof of Lemma~\ref{l.abstract}]
To show $\mathcal E$ consists of pairwise eventually different functions, fix distinct $g_0$ and $g_1$ from $\mathcal E$ and suppose $g_i = \edh(f_i)$ for each $i\in \{0,1\}$.
Clearly we can disregard the set 
\[
N = \{ n\in \nat \setdef g_0 (n)= \ed(f_0)(n) \text{ and } g_1(n) = \ed(f_1)(n) \}
\]
as $g_0$ and $g_1$ can only agree on finitely many such $n$.

If $n \notin N$ then it must be the case that for some $i\in \{0,1\}$, $f_i \in \predT$ and $n \in \ho(f_i)$;
suppose $i=0$ for simplicity. 
By \eqref{adf} we may restrict our attention to $\ho(f_0)\setminus \ho(f_1)$ where
 $g_0$ agrees with $f_0$ and $g_{1}$ agrees with $\ed(f_{1})$.
But $f_0$ and $\ed(f_{1})$ can't agree on an infinite subset of $\ho(f_0) \setminus \ho(f_{1})$ by \eqref{T.1}.

It remains to show maximality. So let $f \colon \nat \to \nat$ be given.
If $f \in \predT$ we have $\edh(f) \res\ho(f) =  f\res\ho(f)$ and $\edh(f) \in \mathcal E$ by definition.

If on the other hand $f \notin \predT$ there is $f' \in {}^\nat\nat$ such that $\ed(f')$ agrees with  $f$ on an infinite subset of $\ho(f)$.  As $\edh(f') \in \mathcal E$ it suffices to show $f \mathbin{=^\infty} \edh(f')$.

If $f'  \notin \predT$ as well this is clear as $\edh(f') = \ed(f')$.
If on the contrary $f'  \in \predT$, we have $f \neq f'$ and so
$\ho(f) \cap \ho(f')$ is finite by \eqref{adf}. So $\edh(f')$ agrees with $\ed(f')$ for all but finitely many points in $\ho(f)$
and hence agrees with $f$ on infinitely many points.
\end{proof}

Note that letting $\mathcal T = \{ f \in {}^\nat \nat \setdef$ $f$ is not $\infty$-coherent on $\domain(f) \}$ and $\ho(f) = \domain(f)$ the requirements of Lemma~\ref{l.abstract} are trivially satisfied; but the resulting $\mathcal E$ will not be Borel (only $\Pi^1_1 \vee \Sigma^1_1$).
On the other hand if $\mathcal T$ is $\Delta^1_1$ and $\ho\colon {}^\nat\nat \to \powerset(\nat)$ is $\Sigma^1_1$, then $\mathcal E$ is clearly $\Sigma^1_1$,
and in fact it follows that $\mathcal E$ is $\Delta^1_1$ in this case because\footnote{In this context, the much more general Theorem 1.4.23 in \cite[p.~15]{miller:intro} deserves mention; compare also \cite[35.10,~p.~285]{kechris}.} it is a \emph{medf} and so 
\[
h \notin \mathcal E \iff (\exists g\in {}^\nat\nat)\; h \neq g \wedge h =^\infty g \wedge g \in \mathcal E.
\]
(Of course the function $\ho\colon {}^\nat\nat \to \powerset(\nat)$ is also automatically $\Delta^1_1$.)
We may view the task at hand to be: find a reasonably effective process producing from a function $f$ either a
subset of $\domain(f)$ where $f$ agrees with some $\ed(f')$ or a set $\ho(f)\subseteq \domain(f)$ on which $f$ can be seen effectively to not be $\infty$-coherent. 

\medskip

From this we can sketch what is arguably the core of Horowitz and Shelah's construction from \cite{medf-borel}.
The present author has not verified whether their  construction yields an arithmetic family.

\begin{proof}[Proof of Theorem~\ref{t.borel}]
Given $f \colon \nat \to \nat$ define a coloring of unordered pairs from $\nat$ as follows (supposing without loss of generality that $k < k'$):
\[
c(\{ k,k'\}) = \begin{cases} 0&\text{ if $lh(s_{f(k)})=k$, $lh(s_{f(k')})=k'$, and $s_{f(k)} \subseteq s_{f(k')}$,} \\
1 &\text{ otherwise.}
\end{cases}
\]
Let $\predT$ be such that for every $f\in \predT$ there is an infinite set $X \subseteq \domain(f)$ which is $1$-homogeneous, i.e.,  $c$ assigns the color $1$ to every unordered pair from $X$, and for every $f\notin \predT$ there is an infinite $0$-homogeneous $X \subseteq \domain(f)$.
Then  \eqref{T.2} holds. 
For $f \in \predT$ let $\ho(f)$ be some infinite $1$-homogeneous $X \subseteq \domain(f)$;  
for $f \notin \predT$ let $\ho(f) = \domain(f)$.
Then \eqref{T.1} and \eqref{adf} hold by definition and by Lemma~\ref{l.abstract}, $\mathcal E$ is a \emph{medf}.

By the proof of the Infinite Ramsey Theorem, the set $\predT$ can be chosen to be $\Delta^1_1$ and the function $\ho \colon {}^\nat\nat \to \powerset(\domain(f))$
can be chosen to be $\Sigma^1_1$. 
Thus $\mathcal E$ as defined in Lemma~\ref{l.abstract} is $\Delta^1_1$.
\end{proof}
In the next section, we essentially replace the appeal to the Infinite Ramsey Theorem by a simple instance of the law of excluded middle.

\section{A maximal eventually different family with a simple definition}\label{s.main}\label{s.construction}

We now give a simpler construction of a family satisfying the requirements of Lemma~\ref{l.abstract}.
\begin{dfn}[The \emph{medf} $\mathcal E$]~
\begin{enumerate}[A.]
\item Let $f\colon \nat\to \nat$. Define a binary relation $\order$ on $\nat$ by
\[
 m \order m' \iff \Big[ \big( \lh(s_{f(m)}) = m \big) \wedge \big(\lh(s_{f(m')}) = m'\big) \wedge \big( s_{f(m)} \subsetneq s_{f(m')}\big) \Big]
\]
\item Let $\predT$ be the set of $f\colon \nat\to\nat$ such that 
\begin{equation}\label{case1}
(\forall n \in \domain(f))(\exists m \in \domain(f) \setminus n) (\forall m' \in \domain(f)\setminus m)\;
\neg( m \order m' ) 
\end{equation}
We also say $f$ \emph{is tangled} to mean $f \in \predT$.
\item For $f\notin \predT$, define $\ho(f)$ to be $\domain(f)$ and for $f \in \predT$ define
\[
\ho(f) = \{ m \in \domain(f) \setdef (\forall m' \in \domain(f)\setminus m)\;
\neg( m \order m' ) \}.
\]

\item Let 
$\mathcal E $ be defined  from $\predT$ and $\ho$ as in Lemma~\ref{l.abstract}, i.e., 
\[
\mathcal E = \{\; \edh(f) \setdef f\in {}^{\nat}\nat \}
\]
where $\edh(f)$ is the function defined as in \eqref{e.e.ddot}:
\begin{equation*}
\edh(f) = \begin{cases}
\edb(f,\ho(f)) & \text{ if $f \in \predT$,}\\
\ed(f)             & \text{ otherwise.}
\end{cases}
\end{equation*}
\end{enumerate}
\end{dfn}
We want to call the following to the readers attention:
\begin{enumerate}[label=(\roman*),ref=\roman*]
\item\label{inst.adf}
$\{ \ho(f) \setdef f \in {}^\nat\nat \}$ is an almost disjoint family (as $\ho(f) \subseteq \domain(f)$ by definition).
\item\label{inst.T.1} When $f$ is tangled, $\ho(f)$ is an infinite set by \eqref{case1} and 
for no $f'\in {}^\nat\nat$ does $f$ agree with $\ed(f')$ on infinitely many (or in fact, just two) points in $\ho(f)$---i.e., $f$ is not $\infty$-coherent on $\ho(f)$.
\end{enumerate}

\begin{lem}\label{main}
The set $\mathcal E$ is a maximal eventually different family.
\end{lem}
\begin{proof}
We show that Lemma~\ref{l.abstract} can be applied. Requirements \eqref{adf} and 
\eqref{T.1} hold by \eqref{inst.adf} and \eqref{inst.T.1} above.
For \eqref{T.2}, suppose $f$ is not tangled, i.e.,
\begin{equation*}\label{case2}
(\exists n \in \domain(f))(\forall m \in \domain(f) \setminus n) (\exists m' \in \domain(f)\setminus m)\;
 m \order m'. 
\end{equation*}
Let $m_0$ be the least witness to the leading existential quantifier above;
by recursion let
$m_{i+1}$ be the least $m'$ in $\domain(f)$ above $m_i$ such that 
 $m_i\order m'$.
Letting 
$f'  = \bigcup \{ s_{f(m_i)} \setdef i \in \nat \}$ yields a well-defined function in ${}^\nat\nat$ such that $f \mathbin{=^\infty} \ed(f')$, i.e., $f$ is $\infty$-coherent on $\ho(f)$.
\end{proof}

It is obvious that $\mathcal E$ is $\Delta^1_1$. 
We now show a stronger result.
\begin{lem}\label{l.pi04}
The set $\mathcal E$ is in the Boolean algebra generated by the $\Sigma^0_3$ sets in ${}^\nat\nat$.
\end{lem}
\begin{proof}
By construction $g \in \mathcal E$ if and only if the following holds of $g$ (see Remark~\ref{r.f.recursive}):
\begin{enumerate}[ref=\Roman*,label=(\Roman*)]
\item $(\forall n \in \nat) \; \lh(s_{g(2n)}) = 2n$, and
\item $(\forall n \in \nat) (\forall m \leq n)\; s_{g(2m)} \subseteq s_{g(2n)}$, and letting $f = \bigcup_{n \in 2\nat} s_{g(2n)}$, 
\item\label{i.inE.tangledcase} either the following three requirements hold: 
                    \begin{enumerate}[label=(\alph*),ref=\ref{i.inE.tangledcase}\alph*]
                    \item\label{i.tangled} $f$ is tangled and
                    \item $(\forall n \in \nat) \; n \in \ho(f) \Rightarrow g(n) = f(n)$ and
                    \item $(\forall n \in \nat) \; n \notin \ho(f) \Rightarrow g(n) = \ed(f)(n)$;
                    \end{enumerate}
\item\label{i.inE.untangledcase} or both of the following hold: 
                    \begin{enumerate}[label=(\alph*),ref=\ref{i.inE.untangledcase}\alph*]
                    \item $f$ is not tangled and
                    \item $(\forall n \in \nat) \; g(n) = \ed(f)(n)$.
                    \end{enumerate}
\end{enumerate} 
As $\ho(f)$ is $\Pi^0_1(f)$ for $f \in \predT$ and \eqref{i.tangled} is $\Pi^0_3(f)$, clearly \eqref{i.inE.tangledcase} is $\Pi^0_3(g,f)$. Likewise \eqref{i.inE.untangledcase} is $\Sigma^0_3(g,f)$.
As $f$ is recursive in $g$, \eqref{i.inE.tangledcase} can be expressed by a $\Pi^0_3(g)$ formula and \eqref{i.inE.untangledcase} can be expressed by a $\Sigma^0_3(g)$ formula (substitute each expression of the form $f(n)=m$ by $s_{g(2n +2)}(n) = m$ and $f \res n$ by $s_{g(2n)}\res n$).
\end{proof}

\section{Mangling away existential quantifiers}\label{s.mangling}

We use the following lemma to reduce the complexity of the family $\mathcal E$.
\begin{lem}
Let $\xi < \omega_1$. Suppose there is a $\Pi^0_{\xi+2}$ maximal eventually different family. Then there is a $\Pi^0_{\xi+1}$  maximal eventually different family.
\end{lem}
\begin{proof}
Suppose
\[
f \in \mathcal E \iff (\forall n \in \nat)(\exists m \in \nat)\; \Psi(n,m,f). 
\]
where $\Psi(n,m,f)$ is $\Pi^0_\xi$.
For each $f \in \mathcal E$ let $g_f \colon \nat \to \nat$ be the function 
such that for each $n\in \nat$, $g_f(n)$ is the least $m$ satisfying $\Psi(n,m,f)$.

We construct a set $\mathcal E^*$ of functions from $\nat$ to $\nat$ as follows.
Given $f \in \mathcal E$, let
$f^* \colon \nat \to \nat$ be the following function: for $n \in \nat$ and $i \in \{0,1\}$ let
\[
f^*(2n+i) = \begin{cases}
f(n) & \text{ for $i=0$;}\\
\#\big(f\res n+1 \conc (g_f \res n+1)\big) & \text{ for $i=1$.}
\end{cases}
\]
It is straightforward to check that $\mathcal E^*$ is a \emph{medf} as every function will agree with an element of $\mathcal E^*$ on infinitely many even numbers.

\medskip

Lastly, $\mathcal E^*$ is $\Pi^0_{\xi+1}$: 
Let $\Psi'(n,m,h)$ denote the formula obtained from $\Psi(n,m,f)$ by replacing each occurrence
of $f(m)=n$ by $h(2m)=n$. 
Clearly $\Psi'$ is $\Pi^0_\xi$.

Let $S_2$ denote the recursive set $\{ m \in \nat\setdef (\exists n\in \nat) \; s_m \in {}^{2n}\nat \}$ and given $m \in S_2$,
write $f_m$ for $s_m\res \frac{\lh(s_m)}{2}$ and $g_m$ for the function $t\colon n \to \nat$ given by $k \mapsto s_n(n+k)$.
In other words, if $m = \#(f\res n+1 \conc g_f \res n+1)$ as above in the definition of $f^*$, then $f_m = f\res n+1$ and $g_m = g_f \res n+1$.
Clearly $m \mapsto f_m$ and $m \mapsto g_m$ are both recursive on $S_2$.

It is straightforward to check that $h \in \mathcal E^*$ if and only if 
for every $n \in \nat$ all of the following hold:
\begin{enumerate}[label=(\roman*),ref=\roman*]
\item\label{i.Estar.pairs}
$h(2n+1) \in S_2 \wedge h(2n) = f_{h(2n+1)}(n)$
\item\label{i.Estar.psi}
$\Psi'( n, g_{h(2n+1)}(n) , h)$
\item\label{i.Estar.g} $(\forall m < g_{h(2n+1)}(n)) \; \neg \Psi'( n, m , h)$.
\end{enumerate}
Requirement \eqref{i.Estar.pairs} is $\Delta^0_1(h)$; \eqref{i.Estar.psi} is $\Pi^0_\xi(h)$ and \eqref{i.Estar.g} is $\Sigma^0_{\xi}(h)$.
So $\mathcal E^*$ is $\Pi^0_{\xi+1}$.
\end{proof}
In fact (but we have no use for this) it is possible to carry out out a similar construction as the above for limit $\xi$. 
This would give a second proof that there is a $\Pi^0_1$ \emph{medf} based on the construction of \emph{any} $\Delta^1_1$ \emph{medf} regardless of its precise complexity, and a version of the above lemma.

\begin{cor}
There is a $\Pi^0_1$ maximal eventually different family.
\end{cor}
\begin{proof}
By Lemma \ref{l.pi04} there is an arithmetic (in fact $\Sigma^0_3\vee\Pi^0_3$) \emph{medf} so we obtain a $\Pi^0_1$ \emph{medf} by the previous lemma. 
\end{proof}

\bibliography{medf-closed}{}
\bibliographystyle{amsplain}

\end{document}

%% file: preamble.tex
\theoremstyle{plain}
\newtheorem{thm}{Theorem}[section]
\newtheorem{lem}[thm]{Lemma}

\newtheorem{cor}[thm]{Corollary}

\newtheorem*{question*}{Question}

\theoremstyle{definition}
\newtheorem{dfn}[thm]{Definition}

\newtheorem{rem}[thm]{Remark}

\theoremstyle{remark}

{\end{minipage}\end{equation}}

\newenvironment{eqpar*}{\begin{equation*}\begin{minipage}{0.8\columnwidth}}%
{\end{minipage}\end{equation*}}

\providecommand{\nat}{\mathbb{N}}

\DeclareMathOperator{\lh}{lh}

\providecommand{\res}{\mathbin{\upharpoonright} }
\providecommand{\conc}{ \mathbin{{}^\frown}}

\providecommand{\Hhier}{\mathbf{H}}

\providecommand{\setdef}{\;|\;}